\documentclass[a4paper]{article}
\title{More two-distance counterexamples to Borsuk's conjecture from strongly regular graphs}
\author{Thomas Jenrich}
\date{2021-06-26}
\begin{document}
\maketitle

\section{Abstract}

In \cite{Bon} (2013) and finally \cite{Bon2} Andriy V. Bondarenko showed
how to construct a two-distance

counterexample to Borsuk's conjecture
from any strongly regular graph whose vertex set is not the union of at
most $f+1$ cliques (sets of pairwise adjacent vertices) where $f$ is the
multiplicity of the second-largest eigenvalue of its adjacency matrix.

He applied that construction to those two graphs that he had been able to
prove to fulfill the condition: From the $G_2(4)$ graph (on 416 vertices)
he got a 65-dimensional two-distance counterexample. From the
$Fi_{23}$ graph (on 31671 vertices) he got a 782-dimensional one and,
by considering certain induced subgraphs, counterexamples in dimensions
781, 780 and 779.

This article presents two other strongly regular graphs fulfilling the
condition, on 28431 and on 2401, resp., vertices. It gives dedicated
counterexamples in dimensions from 781 down to 764 derived from the bigger
graph (that turned out to be an induced subgraph of the $Fi_{23}$ graph)
and a 240-dimensional counterexample derived from the smaller graph.

Several contained propositions rely on the results of (often extensive)
computations, mainly within the computer algebra system GAP. The source
package contains (almost) all used source files.

\section{Introduction}

\subsection{Borsuk's conjecture}

In \cite{Bor} (1933) Karol Borsuk asked whether each bounded set in the
$n$-dimensional Euclidean space can be divided into $n+1$ parts of smaller diameter.
The diameter of a set is defined as the supremum (least upper bound) of the
distances of contained points. Implicitly, the whole set is assumed to
contain at least two points.

The hypothesis that the answer to that question is positive became famous
under the name \emph{Borsuk's conjecture}.

Beginning with Jeff Kahn and Gil Kalai in 1993, several authors have
proved that in certain (almost all) high dimensions such a division is
not generally possible.

\subsection{Preliminaries}

We consider simple loopless finite undirected graphs.

For any graph $\Gamma=(V,E)$, any $a \in V$ and $W \subseteq V$:

$N(\Gamma,a,W) = \{ b \in W : (a,b) \in E \}$ ,
 $n(\Gamma,a,W) = |N(\Gamma,a,W)| $ and

$\Gamma[W]$ is the subgraph of
$\Gamma$ induced by $W$, i.e., the graph whose vertex set is $W$ and whose
edge set is $E \cap (W \times W) $ .

Saying that a vertex $j$ is a non-neighbour of a vertex $i$ means that $i$
and $j$ are neither neighbours nor identical.

For a positive integer $d$ and a point set (vector set) $X$, the
propositions ``X is in dimension d'' and ``X is d-dimensional''
both mean that the exact affine dimension of $X$ is at most $d$.

\subsection{Strongly regular graphs}

A graph $\Gamma=(V,E)$ is called strongly regular with parameter set
$(v,k,\lambda,\mu)$, or shortly a srg$(v,k,\lambda,\mu)$,
iff $|V|=v$ and for all $i, j \in V$
$$ |\{ h \in V : (h,i) \in E \land (h,j) \in E\} | =
\left\{\begin{array}{rl}
k & \mbox{if $i=j$}\\
\lambda & \mbox{if $(i,j) \in E$}\\
\mu & \mbox{otherwise}\end{array}\right.
$$

\subsection{Euclidean representations of strongly regular graphs}

This construction follows \cite{Bon}.

Let $\Gamma=(V,E)$ a srg$(v,k,\lambda,\mu)$ and $A$ its $(0,1)$-adjacency
matrix. Then $A$ has exactly 3 different eigenvalues:
$k$ of multiplicity $1$, the second-largest eigenvalue of multiplicity
$$
f=\frac 12
\left(v-1-\frac{2k+(v-1)(\lambda-\mu)}{\sqrt{(\lambda-\mu)^2+4(k-\mu)}}\right),
$$
and the smallest eigenvalue
$$
s=\frac 12\left(\lambda-\mu-\sqrt{(\lambda-\mu)^2+4(k-\mu)}\right)
$$

In the remaining part of this article we use these notations:
$I$ is the identity matrix of size $v$, $y$ is $A-sI$, $y_i$,
where $i \in V$, are the columns of $y$, and $y_{i,j}$, where
$i,j \in V$, are the entries of $y$, and for any $W \subseteq V$,
the corresponding point set $P(W)$ is $\{ y_i: i \in W \}$.

Remark: Bondarenko used these $y_i$ as intermediate values to derive sets
of $z_i$ and finally $x_i$ by synchronous scaling and moving in order to
get all points onto the unit sphere around the origin but this is not
necessary (here) and could lead to non-integer coordinate values.

The given properties of the eigenvalues imply $\dim P(V) = f$.

For $i, j \in V$
$$
y_{i,j} = \left\{\begin{array}{rl}
-s & \mbox{if $i=j$}\\
1 & \mbox{if $(i,j) \in E$}\\
0 & \mbox{otherwise}\end{array}\right.
$$

For $i\in V$: $y_i$ consists of one -s (at position $i$), $k$ $1$'s, and
$v-k-1$ $0$'s; its norm is the square root of $ s^2 + k $ .

For different $i,j\in V$
$$
\| y_i - y_j \| ^ 2 = \left\{\begin{array}{rl}
 2\times{(k-\lambda-1+(-s-1)^2)} =  2\times{(k-\lambda+s^2+2s)}
 & \mbox{if $(i,j) \in E$}\\
2\times{(k-\mu+s^2)}& \mbox{otherwise}\end{array}\right.
$$

The distance square for the non-adjacent case exceeds the distance square
for the adjacent case by
 $2\times{(\lambda-\mu-2s)} = 2\times{\sqrt{(\lambda-\mu)^2+4(k-\mu)}}$

If the graph is not complete, this excess is positive and we can
conclude:

For any two different $i,j \in V$, the distance of $y_i$ and $y_j$ is
smaller than the diameter of the complete vector set if and only if $i$
and $j$ are neighbours. Thus, for each $W \subseteq V$ the diameter of the
corresponding point set $P(W)$ is smaller than that of $P(V)$ if and only
if $W$ is a clique. And $P(V)$ can be divided into $x$ parts of smaller
diameter if and only if $V$ can be divided into $x$ cliques.

\subsubsection{Euclidean representations for vertex subsets}

We will often consider $P(W)$ for $W \subset V$.
After recognising that $P(V)$ is a counterexample to Borsuk's
conjecture, we will try to find ones in a smaller dimensions.
For that purpose, we will several times prove for subsets
$W_2 \subset W_1$ of $V$, that $ \dim P(W_2) \leq \dim P(W_1)-1$.

This inequality is fulfilled iff there are
$x \in \mathbf{R}^v$ and $ c \in \mathbf{R} $ such that

$\forall i\in W_2 : \langle x, y_i \rangle = c $ and
$\exists i\in W_1 : \langle x, y_i \rangle \neq c $.

As it turned out to be useful in several cases, for any $i\in V$ its
non-neighbourhood

$\{ j \in V : i \neq j \land (i,j) \not\in E \}$
induces a lower-dimensional vector set because for $i,j\in V$
$$
\langle y_i , y_j \rangle = \left\{\begin{array}{rl}
 s^2+k & \mbox{if $i=j$}\\
 \lambda-2s & \mbox{if $(i,j) \in E$}\\
 \mu & \mbox{otherwise}\end{array}\right.
$$

and we know from above that
 $\lambda-2s-\mu = \sqrt{(\lambda-\mu)^2+4(k-\mu)}$ equals 0 only if
the graph is complete (and thus any non-neighbourhood empty).

\subsection{Bondarenko's concrete results}

In \cite{Bon} and \cite{Bon2}, A. Bondarenko considered the
$G_2(4)$ graph, a srg$(416,100,36,20)$, and the $Fi_{23}$ graph, a
srg$(31671,3510,693,351)$, whose corresponding point sets are of
dimensions 65 and 782, resp.

He proved that the sizes of contained cliques cannot exceed 5 and 23,
resp. Because $416/5>83$ and $31671/23=1377$, the corresponding point
sets cannot be divided into less than 84 and 1377, resp., parts of
smaller diameter.

For each of those two graphs, he also gave a construction of a family of
point sets which are dedicated counterexamples for all dimensions that
are larger than that of the respective initial one.

In addition, he considered the vertices of the $Fi_{23}$ graph that are
non-neighbours of 1, 2 and 3, resp., pairwise adjacent fixed vertices.
The corresponding point sets are in dimensions 781, 780 and 779, resp.,
contain (at least) 28160, 25344 and 23040, resp., elements, and thus
cannot be divided into less than 1225, 1102 and 1002, resp., parts of
smaller diameter.
\medskip

Remark:
As the reduction of the decrement of the number of remaining vertices
indicates, the sizes of the common non-neighbourhoods of 4 and 5, resp.,
pairwise adjacent vertices are also large enough to provide counterexamples
to Borsuk's conjecture in dimensions 778 and 777, resp., but for the
concrete numbers (of points and parts) one would have to do actual
countings on the base of actually choosen vertices.

\section{Regular partitions and induced subgraphs}

In the case of the $G_2(4)$ graph, the non-neighbourhood of a fixed vertex
contains exactly 315 vertices. The corresponding point set is in dimension
64 and cannot be divided into less than 63 parts, but this does not prove
the point set to be a counterexample.

But, as shown in \cite{Jen} and \cite{JB}, one can construct a
64-dimensional counterexample from a certain regular partition of the
vertex set of the $G_2(4)$ graph.

As explained in the following, the $Fi_{23}$ graph does have the properties
needed for an analogous construction, finally allowing to derive more than
a dozen different counterexamples in dimensions below of 782.

\subsection{A two-step construction of subsets}

Let $\Gamma=(V,E)$ a strongly regular graph, $A$ its adjacency matrix,
$s$ the smallest eigenvalue of $A$, $y=A-sI$, $W \subseteq V$,
$\{B_1,B_2,B_3,C\}$ a partition of $V$,
$(B_1 \cup B_2) \cap W \ne \emptyset$, $B_3 \cap W \ne \emptyset$,
such that

(1) $\forall g,h \in \{1,2,3\} :
 g \ne h \longrightarrow \forall i \in B_g : n(\Gamma,i,B_h)=0 $

(2) $\forall i \in C : n(\Gamma,i,B_1)=n(\Gamma,i,B_2)=n(\Gamma,i,B_3) $

Recall that for $i, j \in V$
$$
y_{i,j} = \left\{\begin{array}{ll}
-s & \mbox{if $i=j$}\\
1 & \mbox{if $(i,j) \in E$}\\
0 & \mbox{otherwise}\end{array}\right.
$$

For $h \in \{1,2,3\}$ let $x_h$ the vector with index set $V$
whose entries $x_{h,j}$, $j \in V$, are 1 if $j \in B_h$ and 0 otherwise.
Consequently, for any $h \in \{1,2,3\}$ and $i \in V$

$$\langle x_h, y_i \rangle = \sum_{j\in B_h} y_{i,j} =
\left\{\begin{array}{rl}
n(\Gamma,i,B_h)-s > 0 & \mbox{if $j \in B_h$}\\
n(\Gamma,i,B_h) & \mbox{if $j \in C$}\\
n(\Gamma,i,B_h) =0 & \mbox{otherwise}\end{array}\right.
$$

Let $p = x_1 - x_2$.
Combined with (1) and (2) this implies $\forall i \in V$
$$
\langle p, y_i \rangle
= \langle x_1, y_i \rangle - \langle x_2, y_i \rangle
= \\ \left\{\begin{array}{rl}
(n(\Gamma,i,B_1) -s) - 0 = n(\Gamma,i,B_1)-s > 0 & \mbox{if $i \in B_1$}\\
0 - (n(\Gamma,i,B_2)-s) = s-n(\Gamma,i,B_2) < 0 & \mbox{if $i \in B_2$}\\
0 - 0  = 0 & \mbox{if $i \in B_3$}\\
n(\Gamma,i,B_1) - n(\Gamma,i,B_2) = 0 & \mbox{if $i \in C$}\end{array}\right.
$$

We have $\{ i \in W : \langle p, y_i \rangle = 0 \} = (B_3 \cup C) \cap W
= W \setminus (B_1 \cup B_2) \subset W$.

Thus, $ \dim P((B_3 \cup C) \cap W) \leq \dim P(W)-1 $.

Let $q=x_1+x_2 - 2 x_3$.
Combined with (1) and (2) this implies $\forall i \in V$
$$
\langle q, y_i \rangle
= \langle x_1, y_i \rangle + \langle x_2, y_i \rangle - 2 \langle x_3, y_i \rangle
=\\
\left\{\begin{array}{rl}
(n(\Gamma,i,B_1)-s) + 0 -  2 \times 0 = n(\Gamma,i,B_1)-s & \mbox{if $i \in B_1$}\\
0 + (n(\Gamma,i,B_2)-s) -  2 \times 0 = n(\Gamma,i,B_2)-s & \mbox{if $i \in B_2$}\\
0 + 0 - 2 \times (n(\Gamma,i,B_3)-s) = 2\times (s-n(\Gamma,i,B_3)) & \mbox{if $i \in B_3$}\\
n(\Gamma,i,B_1) + n(\Gamma,i,B_2) - 2 \times n(\Gamma,i,B_3) = 0 & \mbox{if $i \in C$}\end{array}\right.
$$

We have $\{ i \in (B_3 \cup C) \cap W : \langle q, y_i \rangle = 0 \} =
 C\cap W \subset  (B_3 \cup C) \cap W $.

Thus, $ \dim P(C \cap W) \leq  \dim P((B_3 \cup C) \cap W) -1 $.

\subsection{Application to the $G_2(4)$ graph}

As discussed in \cite{Jen} and \cite{JB}, one can find (many but
equivalent) partitions $\{B_1,B_2,B_3,C\}$ of the vertex set $V$
with $|B_1|=|B_2|=|B_3|=32$ such that the conditions for the construction are fulfilled in the
case $W=V \land b=20 \land c=8$. The first reduction step results in a
64-dimensional point set of size 352 that cannot be divided into less than
71 parts of smaller diameter. It was the first known 64-dimensional
counterexample to Borsuk's conjecture.
The second reduction step (considered in \cite{Jen}) results
in a 63-dimensional point set of size 320 that can be divided into 64
parts of smaller diameter and therefore is no counterexample.

\subsection{Application to the $Fi_{23}$ graph}

One can find a partition $\{B_1,B_2,B_3,C\}$ of the vertex set $V$ with
$|B_1|=|B_2|=|B_3|=1080$ such that the conditions for the construction are
fulfilled in the case $W=V \land b=351 \land c=120$.

This proof (of this fact) is just an adaption of the first part of the
section 4 of \cite{JB}:

Let ($\Gamma=(V,E)$ the $Fi_{23}$ graph and) $\Sigma$ the $Fi_{24}$ graph, a
srg(306936,31671,3510,3240).
It is well-known (cf. \cite{Hub}) that $\Gamma$ occurs as the local graph
of $\Sigma$. Let $x_0$ and $x_1$ two nonadjacent vertices of
$\Sigma$. We can identify the set of the 31671 neighbours of $x_0$ with $V$;
consequently, the common neighbours of $x_0$ and $x_1$ form a 3240-subset,
say, $B$ of $V$.

The graph $\Sigma$ has a triple cover $3\cdot \Sigma$ on 920808 vertices.
It is distance-transitive with intersection array
$\{31671,28160,2160,1; 1,1080,28160,31671\}$.
We see that $B$ is the disjoint union of three mutually nonadjacent
subsets of size 1080. We call them $B_1$, $B_2$ and $B_3$.
Let $C=V\setminus B$.
According to \cite{JKT}, $3 \cdot \Sigma$ is tight, and (see in particular
Figure A.4 in \cite{JKT}) the partition $\{B1,B2,B3,C\}$ possesses the
required properties.

\subsubsection{(Repeated) application of the two-step subset construction}

As in the just given proof, we take the $Fi_{24}$ graph $\Sigma$ as the
initial object, choose a vertex $x_0$ of $\Sigma$, get the $Fi_{23}$ graph
$\Gamma$ as the subgraph induced by the neighbourhood of $x_0$. We start
with $\Gamma$ and its vertex set $V$ and perform one or more two-step
subset constructions, which we will call rounds and number from 1 onwards.

In round $k$, we choose a vertex $x_k$ of $\Sigma$ that is not adjacent to
$x_0$ and different from $x_0, \dots, x_{k-1}$. As seen above, one can
get the partition $\{B_1,B_2,B_3,C\}$ of $V$ where $C$ is the
non-neighbourhood of $x_k$ in $V$ and $B_1,B_2$ and $B_3$ are the components
of the neighbourhood of $x_k$ in $V$.
We use this partition to derive from a subset $Z_{2k-2}$ of $V$
the subsets $Z_{2k-1}=Z_{2k-2} \cap (B_3 \cup C)$ and
$Z_{2k}=Z_{2k-2} \cap C $ .
Let $Z_0=V$.

Because the numbering of the three components of the neigbourhood of $x_k$
in $V$ is not automatically determined, one has to add the information
which of the components should be denoted as $B_3$ in order to define
$Z_{2k-1}$ completely.

Notice that $\dim P(Z_0) = \dim P(V) = 782$ and
 $$Z_{2k} \neq Z_{2k-1} \neq Z_{2k-2} \Longrightarrow
  \dim P(Z_{2k}) < \dim P(Z_{2k-1}) < \dim P(Z_{2k-2}) . $$

Clearly, the dimension decrements are at least 1.

\medskip

\medskip

\section{The computations}

\subsection{Prehistory}

In subsection 3.3 of \cite{CRS} and of \cite{CRS2} the authors in
particular derived from the group $O(7,3)$, also known as $O_7(3)$, a
srg$(28431,3150,621,315)$, denoted $\Gamma^3_6$,
and a srg$(28431,2880,324,288)$, denoted $\Gamma^3_7$, which seemed to be
the first known SRGs with the respective parameter sets.

The appendix of \cite{CRS2} contains a script for the GAP computer
algebra system, helping the readers to get hands on both graphs. In contrast to the body of \cite{CRS2}, the appendix is
accessible without restriction.

The versions 1 to 4 of this article described a (mainly computational)
partial exploration of $\Gamma^3_6$ and derived counterexamples to Borsuk's
conjecture in dimensions down to 774 from it.

In the survey preprint \cite{BvM}, the authors wrote in particular
(without a proof or a dedicated reference) that $\Gamma^3_6$ is an induced
subgraph of the $Fi_{23}$ graph, more precisely the subgraph of the
$Fi_{24}$ graph induced by the intersection of the neighbourhood of a
random vertex and the non-neighbourhood of a random non-neighbour of that
vertex. In terms from the previous section this means that $\Gamma^3_6$ is
(isomorphic to) $\Gamma[Z_2]$, independent of the choice of $x_0$ and
$x_1$ (implying that counterexamples to Borsuk's conjecture derived from
$\Gamma^3_6$ are also counterexamples derived from the $Fi_{23}$ graph).

The computations considered herein include a successful check of this
assertion.

\subsection{Computation outline and results}

Preliminary note: It is common to speak of the $Fi_{23}$ graph and the $Fi_{24}$ graph,
thereby abstracting from the concrete (labellings of the) vertex sets and
considering isomorphic graphs as equal, even identical. In order to express
that we consider individual graphs and not isomorphism classes of graphs,
we will speak of $Fi_{23}$ graphs and $Fi_{24}$ graphs in the remainder
of this section.

\medskip

From a certain $Fi_{24}$ graph $\Sigma$, we get a $Fi_{23}$ graph $\Gamma$
as the subgraph induced by the neighbourhood of vertex 1. We choose eight
different non-neighbours $e_1,\dots ,e_8$ of vertex 1 and build three
lists from these external (with respect to
$\Gamma$) vertices of $\Sigma$, namely  $L_1=(e_1,e_2,e_3,e_6,e_7)$,

 $L_2=(e_1,e_2,e_3,e_4,e_5,e_6,e_7)$ and
 $L_3=(e_1,e_2,e_3,e_4,e_5,e_6,e_8)$ .

For each of these three lists, we do this:

We take the contained vertices of $\Sigma$ in the given order and perform
the described two-step subset construction rounds. In each performed round
$k$, we number the components of the neighbourhood of the used external
vertex in a way such that $ |Z_{2k-1}| $ is as large as possible, i.e.,

 $ | Z_{2k-2} \cap (B_3 \cup C) | \geq | Z_{2k-2} \cap (B_1 \cup C) |$
 and
 $ | Z_{2k-2} \cap (B_3 \cup C) | \geq | Z_{2k-2} \cap (B_2 \cup C) |$.

We check that $ | Z_{2k} | < | Z_{2k-1}| < |Z_{2k-2}| $. This allows us to
conclude that $ \dim P(Z_i) \leq 782 - i $ for all constructed sets $Z_i$.

Because the first three external vertices in those three lists are the same,
not just the initial $Z_0 = V$ is independent of the used list but also
the subsets $Z_1, \dots ,Z_6$. Their respective computed sizes are
 29511, 28431, 26487, 25515, 23571, 22599.

The sizes of the then following subsets $Z_7$ and so on are:

For $L_1$: 21111, 20367, 19119, 18405 .

For $L_2$: 20979, 20169, 18549, 17739, 16611, 16047, 15087, 14553 .

For $L_3$: 20979, 20169, 18549, 17739, 16611, 16047, 15075, 14589 .

\medskip

As mentioned above, Bondarenko proved in particular that the clique number
(largest clique size) of the $Fi_{23}$ graph is 23.
The results of in some cases rather time-consuming computations imply in
particular: The clique numbers of the subgraphs induced by $Z_2$ and
$Z_6$ are 21 and 19, resp., and for $L_2$ and $L_3$, the clique number
of the subgraph induced by $Z_{10}$ is 18.

These results also establish upper bounds for the clique sizes of the
subgraphs induced by the other constructed subsets.

If $W$ is one of the constructed subsets of $V$, $s$ (an upper bound of)
the clique number of $\Gamma[W]$ and an integer $p < |W| / s$,
then $\Gamma[W]$ cannot be divided into $p$ (or less) cliques.

The following table presents actual values for each of the constructed
subsets of $V$.
In each row, the summarizing value in the last cell is the maximum of
the values of $p$ given in preceding cells in the same row and the
values in the last column in following rows.

\medskip

\begin{tabular}{ | r | l | l | l | l | l | l | l | l |}
\hline
Z &Dim & $L_1$ &                 & $L_2$ &                 & $L_3$& & $>$ \\
\hline
 1&781&$e_1$& $29511/23>1283$ &$e_1$& $29511/23>1283$ &$e_1$& $29511/23>1283$ &1353\\
 2&780&     & $28431/21>1353$ &     & $28431/21>1353$ &     & $28431/21>1353$ &1353\\
 3&779&$e_2$& $26487/21>1261$ &$e_2$& $26487/21>1261$ &$e_2$& $26487/21>1261$ &1261\\
 4&778&     & $25515/21>1214$ &     & $25515/21>1214$ &     & $25515/21>1214$ &1214\\
 5&777&$e_3$& $23571/21>1122$ &$e_3$& $23571/21>1122$ &$e_3$& $23571/21>1122$ &1189\\
 6&776&     & $22599/19>1189$ &     & $22599/19>1189$ &     & $22599/19>1189$ &1189\\
 7&775&$e_6$& $21111/19>1111$ &$e_4$& $20979/19>1104$ &$e_4$& $20979/19>1104$ &1111\\
 8&774&     & $20367/19>1071$ &     & $20169/19>1061$ &     & $20169/19>1061$ &1071\\
 9&773&$e_7$& $19119/19>1006$ &$e_5$& $18549/19>976$  &$e_5$& $18549/19>976$  &1006\\
10&772&     & $18405/19>968$  &     & $17739/18>985$  &     & $17739/18>985$  &985\\
11&771&     &                 &$e_6$& $16611/18>922$  &$e_6$& $16611/18>922$  &922\\
12&770&     &                 &     & $16047/18>891$  &     & $16047/18>891$  &891\\
13&769&     &                 &$e_7$& $15087/18>838$  &$e_8$& $15075/18>837$  &838\\
14&768&     &                 &     & $14553/18>808$  &     & $14589/18>810$  &810\\

\hline
\end{tabular}

\pagebreak

The essence in a sentence:

The Euclidean representations of the subgraphs of $\Gamma$ induced by the
constructed subsets of $V$ establish point sets in dimensions

781, 780, 779, 778, 777, 776, 775, 774, 773, 772, 771, 770, 769 and 768,
resp.,

that cannot be divided into

1353, 1353, 1261, 1214, 1189, 1189, 1111, 1071, 1006, 985, 922, 891, 838
and 810, resp.,

parts of smaller diameter.

\subsection{Computation environment and preconditions}

All computations described herein have been done on a system with
Intel Pentium(R) Dual-Core E5500 at 2.80 GHz and 4 GB RAM, running
(Linux distribution) Lubuntu 20.04 (64 bit)), and with the computer
algebra system GAP (\cite{GAP}, version 4.10.2) and the GAP packages
GRAPE (\cite{Soi}, version 4.8.5) and AtlasRep (\cite{AR}, version 2.1.0)
installed.

Certain older versions of AtlasRep may work as well, but version 1.5.1-2
(installed with GAP 4.10.2) turned out to be too old.
For some tasks, in particular some used in the computation considered here,
GRAPE employs the program Dreadnaut (from the popular graph theoretic
software nauty (by Brendan McKay and Adolfo Piperno, \cite{NAU}).
The installation of GRAPE includes the installation of (an older version
of) Dreadnaut.
For version 4.8.5 of GRAPE, in contrast to older versions, using files
instead of (memory based) strings is default in the transmission
of graph data to Dreadnaut, in order to avoid insufficient memory,
in particular in the case of large graphs as treated here.
Originally, GRAPE 4.8.5 refuses to run on GAP versions before 4.11 .
If you (as me) have GAP version 4.10.2 installed and can't or don't
want to upgrade to a newer version, you can remove that obstacle by
modifying the file \texttt{PackageInfo.g} coming with GRAPE 4.8.5: In the
assignment to the record \texttt{Dependencies}, change the string value
assigned to the field \texttt{GAP} from "\verb+>+=4.11" to "\verb+>+=4.10.2".
As recorded in \cite{GRAPE28}, L. H. Soicher, the author of GRAPE,
replied to my respective question, that he did not see an argument
against this work-around.

\subsection{Input files}

Nine (plain ASCII text) files are expected to be in the working directory
before and during the computations. While GA36.g is from the mentioned
appendix of \cite{CRS2}, the others are included in the source package of
this article.

Script files for the GAP system:

\texttt{GA36OUT.g}, \texttt{Fi23OUT.g}, \texttt{Fi24BOR.g},
\texttt{GA36.g}, \texttt{Fi23.g}, \texttt{Fi24.g} and \texttt{WRIGRA.g}.

(The latter four files are to be called from within one of the first three
files).

Command files for Dreadnaut:

\texttt{Fi24L\_Fi23.DRE} and \texttt{Fi24LR1\_GA36.DRE}

\medskip

Times given in the remainder of this section are CPU core running times,
measured on the computer system described above.

\subsection{Task 1 : GAP processes \texttt{Fi24BOR.g}}

The first action is a call of \texttt{Fi24.g}, in order to construct
a $Fi_{24}$ graph $\Phi_{24}$ from the finite group ${Fi_{24}}'$
and assign it to the variable \texttt{Graf}.

The (labels of) 8 certain different non-neighbours (in the outline denoted
$e_1, \dots, e_8$) of vertex 1 of $\Phi_{24}$ are assigned to the variable
\texttt{Ext}.

Three lists of indices of those external (with respect to the neigbourhood
of vertex 1) vertices are stored in the variable \texttt{elis}, reflecting
the definition of the lists $L_1$, $L_2$ and $L_3$ in the outline.

For each of those three lists, the function \texttt{check\_subset\_counts}
is called, to construct the partitions of the neighbourhood of vertex 1
induced by the external vertices, construct the vertex subsets (in the
outline denoted $Z_1$ and so on), calculate their cardinalities and print
the results (onto the screen, at least in the default case).
(21 seconds)

\medskip

Then an offer to start writing graph data into files and checking clique
size bounds appears. The decision is left to the user because the execution
could take long and needs a considerable amount of memory space (mostly for
GAP itself, much less for the called program Dreadnaut).
(2:06 hours and up to 1700 MB main memory space)

There are seven stages (passes). With the exception of the first one, they
are ordered by increasing system demands (in order to get as much as
possible work done when the system becomes overstressed).

\subsubsection{Stage 1}

Builds the subgraph $G$ of $\Phi_{24}$ induced by the vertex set consisting
of the neighbourhood of vertex 1 and the eight external vertices and
assigns it to the variable \texttt{Graf}. (44 minutes)

That graph contains the complete information that is relevant for the
following stages. Its subgraph induced by the vertices from 1 to 31671
is a $Fi_{23}$ graph. The external vertices $e_1, \dots, e_8$ are mapped
to the other vertices (from 31672 to 31679), keeping their order.

Tests have shown that this action considerably reduces the memory space
demands and the running time of the whole task.

\subsubsection{Stage 2}
Writes the data of the subgraph of $G$ induced by the first 31671 vertices
(isomorphic to the local graph of $\Phi_{24}$) into the file
\texttt{Fi24L.DRE}, using the input format of Dreadnaut.
(31671 vertices, 301134276 bytes)

\subsubsection{Stage 3}

Writes the data of the subgraph of $G$ induced by those of the first 31671
vertices that are non-neighbours of the first external vertex into the file
\texttt{Fi24LR1.DRE}, using the input format of Dreadnaut.
(28431 vertices, 240295550 bytes)

\subsubsection{Stage 4}

Writes the data of the subgraph of $G$ induced by those of the first 31671
vertices that are common non-neighbours of the first 5 external vertices
into the file \texttt{Fi24LR5.A}, using the ASCII version of the DIMACS
graph data file format.
(17739 vertices, 221837476 bytes)

\subsubsection{Stage 5}

Writes the data of the subgraph of $G$ induced by those of the first 31671
vertices that are common non-neighbours of the first 3 external vertices
into the file \texttt{Fi24LR3.A}, using the ASCII version of the DIMACS
graph data file format.
(22599 vertices, 368011696 bytes)

\subsubsection{Stage 6}

Checks that the subgraph of $G$ induced by those of the first 31671
vertices that are common non-neighbours of the first 5 external vertices
does not have a clique of size 19.

\subsubsection{Stage 7}

Checks that the subgraph of $G$ induced by those of the first 31671
vertices that are common non-neighbours of the first 3 external vertices
does not have a clique of size 20.

\subsection{Task 2 : GAP processes \texttt{GA36OUT.g}}

Constructs $\Gamma^3_6$ from the finite group ${O_8}^+(3)$, calculates and
shows the index number of its local graph (at vertex 1) and writes its
data into the file \texttt{GA36.DRE}, using the input format of Dreadnaut.
(239663676 bytes, 2:25 minutes)

\subsection{Task 3 : GAP processes \texttt{Fi23OUT.g}}

Constructs a $Fi_{23}$ graph from the finite group $Fi_{23}$ and writes
its data into the file \texttt{Fi23.DRE}, using the input format of
Dreadnaut.
(300751232 bytes, 3:47 minutes, about 1 GB main memory space)

\subsection{Task 4 : Dreadnaut processes Fi24L\_Fi23.DRE }

Reads the data of two graphs from the files \texttt{Fi24L.DRE}
and \texttt{Fi23.DRE}, resp., and checks that those two graphs are
isomorphic, by canonical labelling of both graphs and checking the
results, assigned to the variables h and h', for identity.
(Duration of reading in and labelling: 96.59 seconds and 58.25 seconds,
resp.)

\subsection{Task 5 : Dreadnaut processes Fi24LR1\_GA36.DRE}

Reads the data of two graphs from the files \texttt{Fi24LR1.DRE}
and \texttt{GA36.DRE}, resp., and checks that those two graphs are
isomorphic, by canonical labelling of both graphs and checking the
results, assigned to the variables h and h', for identity.

(Duration of reading in and labelling: 829.77 seconds seconds and 725.48
seconds, resp., much longer than for the analogous operations in task 4
although the graphs are smaller here)

\subsection{Tasks 7 and 8 (optional) : Independent calculation of clique numbers}

Apply a clique number computing program that accepts input files in
the ASCII version of the DIMACS graph data format to one or both of
\texttt{Fi24LR3.A} and \texttt{Fi24LR5.A} generated in stages 4 and 5,
resp., of task 1. This way, the clique numbers could be found (or verified)
independent of GAP but probably at the cost of (much) more computer
runtime if that program does not make use of symmetries as much as
GAP does.
In the tests, the popular Cliquer software (by Sampo Niskanen and
Patric C. {\"O}sterg{\aa}rd) has been used.

\section {A counterexample from a srg(2401,240,59,20)}

The chapter \emph{Individual graph descriptions} of \cite{BvM}
contains a section on rank 3 graphs on 2401 vertices. One of those
graphs is a srg$(2401,240,59,20)$. These parameters imply that the
multiplicity $f$ of the second-largest eigenvalue is 240. Thus, the
Euclidean representation is 240-dimensional. According to the information
given in that subsection, the largest clique size of that graph is 9.
Because $9 \times 266<2401$, that graph can not be divided into less than
267 cliques, and its Euclidean representation can not be divided into less
than 267 parts of smaller diameter.

Information on a construction of that graph is given in preceding
chapters of \cite{BvM}. The vertices are the $7^4$ vectors of $GF(7)^4$,
and two vertices are adjacent iff their difference is in a certain set of
240 vectors from a projective $(40,4,12,5)$ set of points in the projective
space PG(3,7). An explanation of this type of correspondences can be found
in the popular early survey paper \cite{CK}.

\subsection{Derived counterexamples in lower dimensions?}

The subgraph induced by the non-neighbours of a random vertex contains
exactly $2401-1-240=2160$ vertices. The dimension of its Euclidean
representation is $240-1=239$. It would establish a counterexample to
Borsuk's conjecture if a partition into 240 cliques would be impossible.
Because $2160=240 \times 9$, one can not prove this by a single division.

\subsection{Computations}

The source package of this article contains the Pascal source file
\texttt{SRG2401E.PAS}, written on the base of the construction of the
graph considered here as given in \cite{BvM}.
The projective point set has been pre-calculated and is given as
constant \texttt{p3to7}. The variable \texttt{A} to store the adjacency
relation is initialized during runtime.

For the tests, the Free Pascal compiler has been used, started by entering
the command line

\texttt{fpc -Mtp SRG2401E.PAS}

The vast majority of the running time has been used by an (unsuccessful)
search for a 10-clique and by a check that the graph is indeed strongly
regular with the given parameter set.

The respective durations were about 17 and 51 seconds on the computer
system described above, using compiler version 3.0.4 for x86\_64.
They were about 47 and 170 seconds on a 1 GHz PIII running MS Windows
98SE, using compiler version 2.4.4 for i386.

The search for a 10-clique would take much longer if not the symmetry of
the graph caused by its construction using a difference set would allow to
consider just the cliques containing vertex 1.

The strong regularity does not have to be checked computationally because
it can be easily concluded from the construction and the fact that the
instruction \texttt{check\_intersections(12,5)} has been passed.
But I do not know of a mathematical proof of the fact that that graph does
not have a 10-clique.

\vspace{0.1in}

Author's eMail address: thomas.jenrich@gmx.de

\end{document}